\documentclass[11pt,a4paper]{amsart}
\usepackage[T1]{fontenc}
\usepackage{lmodern}
\usepackage{amsmath,amssymb,amsthm}
\usepackage[margin=25mm]{geometry}
\usepackage{microtype}
\usepackage[hidelinks]{hyperref}

\newcommand{\F}{\mathbb F}
\newcommand{\PG}{\mathrm{PG}}
\newtheorem{theorem}{Theorem}
\newtheorem{corollary}[theorem]{Corollary}
\numberwithin{equation}{section}
\hypersetup{
  pdftitle={Functions determining few directions: the missing case},
  pdfsubject={Functions over finite fields determining fewer than (q+3)/2 directions are linear over a subfield},
  pdfkeywords={finite field, direction, function graph, additive function, Redei polynomial, blocking set}
}

\title[Functions determining few directions]
  {Functions determining few directions:\newline the final missing case}
\author{Peter Sziklai}
\date{2 October, 2026}

\begin{document}

\begin{abstract}
By results of Blokhuis, Ball, Brouwer, Storme and Sz\H{o}nyi, and of
Ball, if a function $f: \F_q\to\F_q$ determines fewer than $(q+3)/2$
directions, then $f-f(0)$ is linear over a subfield of $\F_q$, except
possibly when $q$ is even and $f$ determines exactly $q/2+1$
directions. We settle this case: $f-f(0)$ is then additive. Hence every
small minimal blocking set of R\'edei type in $\PG(2,q)$ is linear.
\end{abstract}

\maketitle

\section{Introduction}

Let $q=p^n$ with $p$ prime. For $f: \F_q\to\F_q$ let
$U_f=\{(x,f(x)):x\in\F_q\}$ be its graph, let
$D_f=\{(f(x)-f(z))/(x-z):x\ne z\}$ be the set of directions (slopes)
it determines, and put $N=|D_f|$. Let $s=p^e$ be the largest power of
$p$ such that every line with slope in $D_f$ meets $U_f$ in $0$ or in a
multiple of $s$ points.

\begin{theorem}\label{thm:complete}
If $N<(q+3)/2$, then $s\ge2$, $\F_s$ is a subfield of $\F_q$, the map
$x\mapsto f(x)-f(0)$ is $\F_s$-linear, and either $s=q$ and $N=1$, or
$q/s+1\le N\le(q-1)/(s-1)$.
\end{theorem}

Theorem~\ref{thm:complete} is due to Blokhuis, Ball, Brouwer, Storme
and Sz\H{o}nyi~\cite{BBBSS} and Ball~\cite{Ball}, except in one case.
By \cite[Theorem~1.1]{Ball}, either $s=1$ and $N\ge(q+3)/2$, or $s=q$
and $N=1$, or $\F_s$ is a subfield of $\F_q$ and
$q/s+1\le N\le(q-1)/(s-1)$; moreover, $f-f(0)$ is $\F_s$-linear
whenever $s>2$. (The cases $s=2,3$, left open in~\cite{BBBSS}, were
completed in~\cite{Ball}; for $s=2$ this meant proving $N\ge q/2+1$.)
So if $N<(q+3)/2$, then $s\ge2$, and only $s=2<q$ remains (if $s=q$,
then $N=1$ and $f$ is affine). Then $p=2$,
and $q/2+1\le N<(q+3)/2$ forces $N=q/2+1$. (The same bounds show that
for $q=2^n\ge4$ no function determines exactly $q/2$ directions:
$s=1,2$ are excluded, and $s\ge4$ gives $N\le(q-1)/3$.) This case was
raised by Csajb\'ok~\cite[Problem~1]{Csajbok}, and settling it proves
Theorem~\ref{thm:complete}.

\begin{theorem}\label{thm:main}
Let $q=2^n$ with $n\ge2$, and let $f: \F_q\to\F_q$ determine exactly
$q/2+1$ directions. Then $f-f(0)$ is additive; equivalently, $U_f$ is a
translate of an $\F_2$-subspace of $\F_q^2$.
\end{theorem}

The proof below is self-contained. It uses the R\'edei-polynomial
method of~\cite{Redei,BBBSS,Ball}, with ordinary derivatives and a
count of double roots.

\section{Proof of Theorem~\ref{thm:main}}

Throughout, $q=2^n$, so signs play no role. A prime denotes the
derivative with respect to $X$, and $[X^j]F$ is the coefficient of
$X^j$ in $F$. Replacing $f$ by $f-f(0)$ translates $U_f$ and does not
change $D_f$, so we may assume $f(0)=0$. Write $D=D_f$, so
$N=|D|=q/2+1$, and consider the R\'edei polynomial
\[ R(X,Y)=\prod_{x\in\F_q}\bigl(X+xY+f(x)\bigr)=\sum_{j=0}^{q}a_j(Y)X^j,
 \qquad a_j\in\F_q[Y].
\]
For $y,c\in\F_q$ let $m_c(y)$ be the number of $x$ with $f(x)=xy+c$,
that is, the number of points of $U_f$ on the line of slope $y$ through
$(0,c)$. Then
\begin{equation}\label{eq:fiber}
 R(X,y)=\textstyle
 \prod_{c\in\F_q}(X+c)^{m_c(y)}.
\end{equation}
We note four simple facts.
(a)~$a_q=1$, and $a_0=0$ because the factor with $x=0$ is $X$.
(b)~$a_1(Y)=\prod_{x\ne0}(xY+f(x))$, as only one product of $q-1$
constant terms avoids the factor $X$. So the multiplicity of $y$ as a
root of $a_1$ is $m_0(y)-1$, the number of $x\ne0$ with $f(x)=xy$.
(c)~If $y\notin D$, then all $m_c(y)\le1$, so $R(X,y)=X^q+X$ and
$a_j(y)=0$ for $2\le j\le q-1$.
(d)~The homogeneous part of $R$ of degree $q$ is
$\prod_x(X+xY)=X^q+XY^{q-1}$, so $\deg a_j\le q-j-1$ for $2\le j\le q-1$.

\smallskip\noindent\emph{Step 1: a degree gap.}
By (c), each $a_j$ with $2\le j\le q-1$ vanishes on $\F_q\setminus D$,
a set of $q-N=q/2-1$ elements. If $j>q/2$, then $\deg a_j\le q/2-2$
by (d), so $a_j=0$. With (a), this gives
\begin{equation}\label{eq:gap}
 R(X,Y)=X^q+\textstyle\sum_{j=1}^{q/2}a_j(Y)X^j.
\end{equation}

\smallskip\noindent\emph{Step 2: even intersections.}
Fix $y\in D$ and write $P(X)=R(X,y)=X^q+g(X)$, where $\deg g\le q/2$ by
\eqref{eq:gap}. By \eqref{eq:fiber}, $P'=g'$ is divisible by
$\prod_c(X+c)^{m_c(y)-1}$, the product taken over the $c$ with
$m_c(y)\ge1$. As $X^q+X=\prod_{c\in\F_q}(X+c)$, we get
$P\mid(X^q+X)P'$; since $X^q+X=P+(X+g)$, also $P\mid(X+g)P'$. The
latter product has degree at most $q/2+(q/2-1)<q$, so it is zero. If
$P'\ne0$, then $g=X$ and $P=X^q+X$, so all $m_c(y)\le1$, contradicting
$y\in D$. Thus $P'=0$: $P$ is a square, and every $m_c(y)$ is even.
\emph{Every line with slope in $D$ meets $U_f$ in an even number of
points}; in particular $s\ge2$.

Consequently, for odd $j\ge3$, the polynomial $a_j$ vanishes on $D$
(as $[X^j]P=0$) and on $\F_q\setminus D$ (by (c)); since
$\deg a_j<q$, we get $a_j=0$. Also $a_1(y)=[X]P=0$ for $y\in D$. Hence
\begin{equation}\label{eq:even}
 R(X,Y)=X^q+G(X,Y)+a_1(Y)X,\qquad
 G=\textstyle\sum_{j}a_j(Y)X^j\ \ (j\ \text{even},\ 2\le j\le q/2).
\end{equation}

\smallskip\noindent\emph{Step 3: two-point lines through the origin.}
Let $E$ be the set of $y\in D$ with $m_0(y)=2$, i.e.\ the line of slope
$y$ through $(0,0)$ contains exactly one further point of $U_f$. Each
point $(x,f(x))$ with $x\ne0$ lies on the line of slope $f(x)/x\in D$
through the origin, so $\sum_{y\in D}(m_0(y)-1)=q-1$. By Step~2,
$m_0(y)\ge2$ for every $y\in D$, and $m_0(y)\ge4$ if $y\notin E$. Hence
\begin{equation}\label{eq:E}
 q-1\ge|E|+3(N-|E|),\qquad\text{so}\qquad
 |E|\ge(3N-q+1)/2=q/4+2.
\end{equation}
By (b), every $y\in E$ is a simple root of $a_1$, i.e.\ $a_1'(y)\ne0$.

\smallskip\noindent\emph{Step 4: an identity at the directions in $E$.}
Fix $y\in E$. As $a_1(y)=0$, \eqref{eq:even} gives
$P(X)=R(X,y)=X^q+g(X)$ with $g(X)=G(X,y)$. Let
$H(X)=\frac{\partial R}{\partial Y}(X,y)
=\sum_{x}x\prod_{z\ne x}(X+zy+f(z))$. By \eqref{eq:fiber}, every
summand, hence $H$, is divisible by $\prod_c(X+c)^{m_c(y)-1}$.
Moreover, each root $c$ of $P$ lies in $\F_q$, so $0=c^q+g(c)=c+g(c)$,
and $X+c$ divides $X+g$. Hence $P\mid(X+g)H$. By \eqref{eq:gap},
$H=\sum_{j\le q/2}a_j'(y)X^j$ and $X+g$ both have degree at most $q/2$,
so
\begin{equation}\label{eq:scalar}
 (X+g)H=bP\qquad\text{for some }b\in\F_q.
\end{equation}
By \eqref{eq:even}, $g'=0$, and $H'=a_1'(y)=:\gamma\ne0$, since
$\partial G/\partial Y$ has only even $X$-exponents. Differentiating
\eqref{eq:scalar}, and using $P'=0$ and $(X+g)'=1$, we get
$H+\gamma(X+g)=0$, i.e.\ $H=\gamma(X+g)$. Substituting back,
\begin{equation}\label{eq:key}
 bP=\gamma(X+g)^2=\gamma X^2+\gamma\,g(X)^2 .
\end{equation}

\smallskip\noindent\emph{Step 5: cont'd, power-of-2 exponents.}
Write $g=\sum_j g_jX^j$; by \eqref{eq:even}, $g_j=0$ unless $j$ is even
and $2\le j\le q/2$, and $g(X)^2=\sum_j g_j^2X^{2j}$. Suppose that
$g_j\ne0$ for some $j$ that is not a power of two, and take the largest
such $j$. Then $2j$ is not a power of two, so $2j\notin\{2,q\}$, and
$g_{2j}=0$ by the choice of $j$. Comparing the coefficients of
$X^{2j}$ in \eqref{eq:key} gives $0=\gamma g_j^2$, a contradiction.
Hence $g$, and with it $H=\gamma(X+g)$, contains only monomials
$X^{2^i}$.

\smallskip\noindent\emph{Step 6: double roots.}
Let $2\le j\le q/2$, where $j$ is not a power of two. For $y\in E$,
Step~5 gives $a_j(y)=[X^j]P=0$ and $a_j'(y)=[X^j]H=0$, so $(Y+y)^2$
divides $a_j$. By (c), $a_j$ also vanishes on $\F_q\setminus D$, which
is disjoint from $E$. By \eqref{eq:E}, the number of roots of $a_j$,
counted with multiplicity, is at least
$2|E|+(q-N)\ge(q/2+4)+(q/2-1)=q+3>\deg a_j$. Hence $a_j=0$, and
$R(X,Y)=X^q+\sum_{i=0}^{n-1}c_i(Y)X^{2^i}$ with $c_i\in\F_q[Y]$ is an
additive polynomial in $X$.

\smallskip\noindent\emph{Step 7: additivity of $f$.}
In $\F_q[Y]$, the polynomial $R$ has the $q$ roots $r_x=xY+f(x)$,
$x\in\F_q$. They are pairwise distinct, as their coefficients of $Y$
differ, so they are all the roots of $R$. Since $R$ is additive in $X$,
the sum of two roots is a root: for $x,z\in\F_q$, $r_x+r_z=r_w$ for
some $w$. Comparing the coefficients of $Y$ gives $w=x+z$, and
comparing the constant terms gives $f(x)+f(z)=f(x+z)$. \qed

\smallskip
The hypothesis $N=q/2+1$ is used in Step~1 and in the counts of Steps~3
and~6. For smaller $N$ the bound on $|E|$ is too weak, but those cases
are covered by~\cite{Ball}.

\section{R\'edei-type blocking sets}

A blocking set $B$ in $\PG(2,q)$ is \emph{small} if $|B|<3(q+1)/2$,
and of \emph{R\'edei type} if some line $\ell$ meets it in $|B|-q$
points. Let $B$ be small, minimal and of R\'edei type. Take $\ell$ as
the line at infinity, and a point of $\ell\setminus B$ (it exists, as
$|B\cap\ell|<q+1$) as the vertical direction. The $q$ affine lines
through this point each contain a point of $B$, and they partition the
$q$ points of $B\setminus\ell$; so $B\setminus\ell=U_f$ for some $f$.
A point of $\ell$ with slope in $D_f$ lies on an affine line through
two points of $U_f$, so another affine line through it misses $U_f$,
and the point must be in $B$; by minimality, $B\cap\ell=D_f$. Thus
$B=U_f\cup D_f$, $|B|=q+N$, and $B$ is small exactly when
$N<(q+3)/2$. After a translation, $f(0)=0$. If $f$ is $\F_s$-linear,
then $D_f=\{f(x)/x:x\ne0\}$, and
$B=\{\langle(x,f(x),z)\rangle:x\in\F_q,\ z\in\F_s,\ (x,z)\ne(0,0)\}$
is an $\F_s$-linear set. So Theorem~\ref{thm:complete} gives the
following.

\begin{corollary}\label{cor:redei}
Every small minimal blocking set of R\'edei type in $\PG(2,q)$ is an
$\F_s$-linear set for some subfield $\F_s$ of $\F_q$ with $s\ge2$.
\end{corollary}

This confirms the Linearity Conjecture~\cite{Sziklai} for blocking
sets of R\'edei type; the only new case is $q$ even and $|B|=3q/2+1$.

\bigskip

\noindent{\bf Final remark.} 
AI tools were used to proofread and refine the text of this paper only, and the author takes full responsibility for the results. However, like many researchers, he has been working on related problems using AI extensively, which may have influenced his ideas.

\end{document}